\documentclass[12pt,leqno]{article}
\usepackage{amsmath}
\usepackage{amssymb}
\usepackage{latexsym}
\usepackage{amsfonts}
\usepackage{amsmath}
\usepackage{amsthm}
\usepackage{epsfig}
\usepackage{setspace}
\numberwithin{equation}{section}
\numberwithin{figure}{section}
\newtheorem{theorem}{Theorem}
\newtheorem{prop}[theorem]{Proposition}
\numberwithin{theorem}{section}
\newtheorem{lemma}[theorem]{Lemma}
\newtheorem{corollary}[theorem]{Corollary}
\newtheorem{defin}[theorem]{Definition}
\newtheorem{conjecture}[theorem]{Conjecture}

\newcommand{\E}{{\bf E}}

\newcommand{\fatp}{{\bf P}}

\newcommand{\Hyp}{{\mathbb H}^2}
\newcommand{\Hypn}{{\mathbb H}^n}

\begin{document}
\title{The number of unbounded components in the Poisson Boolean model of continuum percolation in hyperbolic space}

\author{Johan Tykesson\thanks{Department of Mathematical Sciences,
    Division of Mathematical Statistics, Chalmers University of
    Technology, S-41296 G\"oteborg, Sweden. E-mail: {\tt
      johant@math.chalmers.se}. Research supported by the Swedish
    Natural Science Research Council.}}

\maketitle

\begin{abstract}
  We consider the Poisson Boolean model of continuum percolation with
  balls of fixed radius $R$ in $n$-dimensional hyperbolic space
  $\Hypn$. Let $\lambda$ be the intensity of the underlying Poisson
  process, and let $N_C$ denote the number of unbounded components in
  the covered region. For the model in any dimension we show that
  there are intensities such that $N_C=\infty$ a.s. if $R$ is big
  enough. In $\Hyp$ we show a stronger result: for any $R$ there are
  two intensities $\lambda_c$ and $\lambda_u$ where
  $0<\lambda_c<\lambda_u<\infty$, such that $N_C=0$ for $\lambda\in
  [0,\lambda_c]$, $N_C=\infty$ for $\lambda\in(\lambda_c,\lambda_u)$
  and $N_C=1$ for $\lambda\in [\lambda_u,\infty)$.
\end{abstract}

\noindent

\noindent{\bf Keywords and phrases:\/} continuum percolation, phase
transitions, hyperbolic space

\noindent{\bf Subject classification:\/} 82B21, 82B43

\section{Introduction}

We begin by describing the fixed radius version of the so called
\emph{Poisson Boolean} model in ${\mathbb R}^n$, arguably the most
studied continuum percolation model. For a detailed study of this
model, we refer to \cite{meester}. Let $X$ be a Poisson point process
in ${\mathbb R}^n$ with some intensity $\lambda$. At each point of
$X$, place a closed ball of radius $R$. Let $C$ be the union of all
balls, and $V$ be the complement of $C$. The sets $V$ and $C$ will be
referred to as the \emph{vacant} and \emph{covered} regions. We say
that \emph{percolation occurs} in $C$ (respectively in $V$) if $C$
(respectively $V$) contains unbounded (connected) components. For the
Poisson Boolean model in ${\mathbb R}^n$, it is known that there is a
\emph{critical intensity} $\lambda_c\in(0,\infty)$ such that for
$\lambda<\lambda_c$, percolation does not occur in $C$, and for
$\lambda>\lambda_c$, percolation occurs in $C$. Also, there is a
critical intensity $\lambda_c^{*}\in(0,\infty)$ such that percolation
occurs in $V$ if $\lambda<\lambda_c^{*}$ and percolation does not
occur if $\lambda>\lambda_c^{*}$. Furthermore, if we denote by $N_C$
and $N_V$ the number of unbounded components of $C$ and $V$
respectively, then it is the case that $N_C$ and $N_V$ are both almost
sure constants which are either $0$ or $1$. In ${\mathbb R}^2$ it is
also known that $\lambda_c=\lambda_c^{*}$ and that at $\lambda_c$,
percolation does not occur in $C$ or $V$. For $n\ge 3$, Sarkar
\cite{sarkar} showed that $\lambda_c<\lambda_c^{*}$, so that there
exists an interval of intensities for which there is an unbounded
component in both $C$ and $V$.

It is possible to consider the Poisson Boolean model in more exotic
spaces than ${\mathbb R}^n$, and one might ask if there are spaces for
which several unbounded components coexist with positive probability.
The main results of this paper is that this is indeed the case for
$n$-dimensional hyperbolic space $\Hypn$. We show that there are
intensities for which there are almost surely infinitely many
unbounded components in the covered region if $R$ is big enough. In
$\Hyp$ we also show the existence of three distinct phases regarding
the number of unbounded components, for any $R$. It turns out that the
main difference between ${\mathbb R}^n$ and $\Hypn$ which causes this,
is the fact that there is a linear isoperimetric inequality in
$\Hypn$, which is a consequence of the constant negative curvature of
the spaces.  In $\Hyp$, the linear isoperimetric inequality says that
the circumference of a bounded simply connected set is always bigger
than the area of the set.

The main result in $\Hyp$ is inspired by a theorem due to Benjamini
and Schramm. In \cite{itai1} they show that for a large class of
nonamenable planar transitive graphs, there are infinitely many
infinite clusters for some parameters in Bernoulli bond percolation.
For $\Hyp$ we also show that the model does not percolate on
$\lambda_c$. The discrete analogue of this theorem is due to
Benjamini, Lyons, Peres and Schramm and can be found in
\cite{itai2}. It turns out that several techniques from the
aforementioned papers are possible to adopt to the continuous
setting in $\Hyp$.

There is also a discrete analogue to the main result in $\Hypn$. In
\cite{paksmirnova}, Pak and Smirnova show that for certain Cayley
graphs, there is a non-uniqueness phase for the number of unbounded
components. In this case, while it is still possible to adopt their
main idea to the continuous setting, it is more difficult than for
$\Hyp$.

The rest of the paper is organized as follows. In section
\ref{discrete} we give a very short review of uniqueness and
non-uniqueness results for infinite clusters in Bernoulli
percolation on graphs (for a more extensive review, see the survey
paper \cite{jonasson2}), including the results by Benjamini, Lyons,
Peres, Schramm, Pak and Smirnova. In section \ref{hypdisc} we review
some elementary properties of $\Hypn$. In section \ref{poihyp} we
introduce the model, and give some basic results. Section
\ref{nochyp} is devoted to the proof of the main result in $\Hyp$
and section \ref{nochypn} is devoted to the proof of the main
theorem for the model in $\Hypn$.

\section{Non-uniqueness in discrete percolation} \label{discrete}

Let $G=(V,E)$ be an infinite connected transitive graph with vertex
set $V$ and edge set $E$. In $p$-Bernoulli bond percolation on $G$,
each edge in $E$ is kept with probability $p$ and deleted with
probability $1-p$, independently of all other edges. All vertices are
kept. Let $\fatp_p$ be the probability measure on the subgraphs of $G$
corresponding to $p$-Bernoulli percolation. (It is also possible to
consider $p$-Bernoulli site percolation in which it is the vertices
that are kept or deleted, and all results we present in this section
are valid in this case too.)  In this section, $\omega$ will denote a
random subgraph of $G$.  Connected components of $\omega$ will be
called \emph{clusters}.

Let $I$ be the event that $p$-Bernoulli bond percolation contains
infinite clusters.  One of the most basic facts in the theory of
discrete percolation is that there is a critical probability
$p_c=p_c(G)\in [0,1]$ such that ${\bf P}_p(I)=0$ for $p<p_c(G)$ and
${\bf P}_p(I)=1$ for $p>p_c(G)$. What happens on $p_c$ depends on
the graph. Above $p_c$ it is known that there is $1$ or $\infty$
infinite clusters for transitive graphs. If we let $p_u=p_u(G)$ be
the infimum of the set of $p\in[0,1]$ such that $p$-Bernoulli bond
percolation has a unique infinite cluster, Schonmann
\cite{schonmann} showed for all transitive graphs, one has
uniqueness for all $p>p_u$. Thus there are at most three phases for
$p\in[0,1]$ regarding the number of infinite clusters, namely one
for which this number is $0$, one where the number is $\infty$ and
finally one where uniqueness holds.

A problem which in recent years has attracted much interest is to
decide for which graphs $p_c<p_u$. It turns out that whether a graph
is \emph{amenable} or not is central in settling this question:

For $K\subset V$, the {\it inner vertex boundary} of $K$ is defined as
$\partial_{V}K:=\{y\in K:\exists x\notin K,\,[x,y]\in E\}$. The {\it
  vertex-isoperimetric} constant for $G$ is defined as
$\kappa_V(G):=\inf_{W}\frac{|\partial_V W|}{|W|}$ where the infimum
ranges over all finite connected subsets $W$ of $V$. A bounded degree
graph $G=(V,E)$ is said to be \emph{amenable} if $\kappa_V(G)=0$.

Benjamini and Schramm \cite{itai4} have made the following general
conjecture:

\begin{conjecture}
\label{schrammconjecture}
If $G$ is transitive, then $p_u>p_c$ if and only if $G$ is
nonamenable.
\end{conjecture}

Of course, one direction of the conjecture is the well-known theorem
by Burton and Keane \cite{burtonkeane} which says that any transitive,
amenable graph $G$ has a unique infinite cluster for all $p>p_c$.

The other direction of Conjecture \ref{schrammconjecture} has only
been partially solved. Here is one such result that will be of
particular interest to us, due to Benjamini and Schramm \cite{itai1}.
This can be considered as the discrete analogue to our main theorem in
$\Hyp$.  First, another definition is needed.

\begin{defin}
\label{oneend}
Let $G=(V,E)$ be an infinite connected graph and for $W\subset V$ let $N_W$
be the number of infinite clusters of $G\setminus W$. The
number $\sup_W N_W$ where the supremum is taken over all finite $W$ is
called the \emph{number of ends} of $G$.
\end{defin}

\begin{theorem}
\label{blps}
Let $G$ be a nonamenable, planar transitive graph with one end.  Then
$0<p_c(G)<p_u(G)<1$ for Bernoulli bond percolation on $G$.
\end{theorem}

Such a general result is not yet available for non-planar graphs.
However, below we present a theorem by Pak and Smirnova
\cite{paksmirnova} which proves non-uniqueness for a certain class of
Cayley graphs.

\begin{defin}
\label{cayleydef}
Let $\Gamma$ be a finitely generated group and let $S=\{g_{1}^{\pm
  1},...,g_{n}^{\pm 1}\}$ be a finite symmetric set of generators for
$\Gamma$. The (right) \emph{Cayley graph} $\Gamma=\Gamma(G,S)$ is the
graph with vertex set $\Gamma$ and $[g,h]$ is an edge in $ \Gamma$ if
and only if $g^{-1}h\in S$.
\end{defin}

Let $S^k$ be the multiset of elements of $\Gamma$ of the type $g_1
g_2...g_k$, $g_1,...,g_k\in S$ and each such element taken with
multiplicity equal to the number of ways to write it in this way. Then
$S^k$ generates $G$.

\begin{theorem}
\label{paksmirnovasats}
  Suppose $\Gamma=\Gamma(G,S)$ is a nonamenable Cayley-graph and let
  $\Gamma_k=\Gamma(G,S^k)$. Then for $k$ large enough,
  \[p_c(\Gamma_k)<p_u(\Gamma_k).\]
\end{theorem}

Theorem \ref{paksmirnovasats} is the inspiration for our main result
in $\Hypn$.

\section{Hyperbolic space}
\label{hypdisc}
We consider the unit ball model of $n$-dimensional hyperbolic space
${\mathbb H}^n$, that is we consider ${\mathbb H}^n$ as the open
unit ball in ${\mathbb R}^n$ equipped with the hyperbolic metric.
The hyperbolic metric is the metric which to a curve
$\gamma=\{\gamma(t)\}_{t=0}^{1}$ assigns length
\[L(\gamma)=2\int_{0}^{1}\frac{|\gamma'(t)|}{1-|\gamma(t)|^2}dt,\] and
to a set $E$ assigns volume \[\mu(E)=2^n \int_E
\frac{dx_1...dx_n}{(1-|x|^2)^2}.\]

The \emph{linear isoperimetric inequality} for $\Hyp$ says that for
all measurable $A\subset \Hyp$ with $L(\partial A)$ and $\mu(A)$ well
defined,
\begin{equation}
\label{isoperimetricinequality}
\frac{L(\partial A)}{\mu(A)}\ge 1.
\end{equation}

Denote by $d(x,y)$ the hyperbolic distance between the points $x$
and $y$. Let $S(x,r):=\{y\,:\,d(x,y)\le r\}$ be the closed
hyperbolic ball of radius $r$ centered at $x$. In what follows, area
(resp. length) will always mean hyperbolic area (resp. hyperbolic
length). The volume of a ball is given by
\begin{equation}
\mu(S(0,r))=B(n)\int_{0}^r\sinh(t)^{n-1}\,dt
\end{equation}
\noindent where $B(n)>0$ is a constant depending only on the
dimension. We will make use of the fact that for any $\epsilon>0$ there is a constant $K(\epsilon,n)>0$ independent of $r$ such
that
\begin{equation}
\label{iso2} \mu(S(0,r)\setminus S(0,r-\epsilon))\ge K(\epsilon,n)
\mu(S(0,r))
\end{equation}
\noindent for all $r$. For more facts about $\Hypn$, we refer to
\cite{ratcliffe}.

\subsection{Mass transport}

Next, we present an essential ingredient to our proofs in $\Hyp$,
the mass transport principle which is due to Benjamini and Schramm
\cite{itai1}. We denote the group of isometries of $\Hyp$ by
Isom($\Hyp$).

\begin{defin}
\label{diagonalinv}
A measure $\nu$ on $\Hyp\times\Hyp$ is said to be \emph{diagonally
  invariant} if for all measurable $A,\,B\subset\Hyp$ and
$g\in$\emph{Isom($\Hyp$)} \[\nu(g A\times g B)=\nu(A\times
B).\]\end{defin}

\begin{theorem}({\sc Mass Transport Principle in $\Hyp$})
\label{masstransport}
If $\nu$ is a positive diagonally invariant measure on ${\mathbb
  H}^2\times{\mathbb H}^2$ such that $\nu(A\times {\mathbb
  H}^2)<\infty$ for some open $A\subset {\mathbb H}^2$, then
\[\nu(B\times\Hyp)=\nu(\Hyp\times B)\] for all measurable $B\subset
\Hyp$.
\end{theorem}
The intuition behind the mass transport principle can be described as
follows. One may think of $\nu(A\times B)$ as the amount of mass that
goes from $A$ to $B$. Thus the mass transport principle says that the
amount of mass that goes out of $A$ equals the mass that goes into
$A$.

\section{The Poisson Boolean model in hyperbolic space}
\label{poihyp}

\begin{defin}\label{poisson}
  A point process $X$ on $\Hypn$ distributed according to the
  probability measure $\fatp$ such that for $k\in {\mathbb N}$,
  $\lambda\ge 0$, and every measurable $A\subset\Hypn$ one has
\[\fatp[|X(A)|=k]=\mbox{e}^{-\lambda \mu(A)}\frac{(\lambda \mu(A))^{k}}{k!}\]
is called a \emph{Poisson process} with intensity $\lambda$ on $\Hypn$.
Here $X(A)=X\cap A$ and $|\cdot|$ denotes cardinality.
\end{defin}

In the \emph{Poisson Boolean model} in $\Hypn$, at every point of a
Poisson process $X$ we place a ball with fixed radius $R$. More
precisely, we let $C=\bigcup_{x\in X}S(x,R)$ and $V=C^c$ and refer to
$C$ and $V$ as the covered and vacant regions of $\Hypn$ respectively.
For $A\subset\Hypn$ we let $C[A]:=\bigcup_{x\in X(A)}S(x,R)$ and
$V[A]:=C[A]^c$. For $x,\,y\in\Hypn$ we write $x\leftrightarrow y$ if
there is some curve connecting $x$ to $y$ which is completely covered
by $C$. Let $d_C(x,y)$ be the length of the shortest curve connecting
$x$ and $y$ lying completely in $C$ if there exists such a curve,
otherwise let $d_C(x,y)=\infty$. Similarly, let $d_V(x,y)$ be the
length of the shortest curve connecting $x$ and $y$ lying completely
in $V$ if there is such a curve, otherwise let $d_V(x,y)=\infty$. The
collection of all components of $C$ is denoted by ${\mathcal C}$ and
the collection of all components of $V$ is denoted by ${\mathcal V}$.
Let $N_C$ denote the number of unbounded components in $C$ and $N_V$
denote the number of unbounded components in $V$. Next we introduce
four \emph{critical intensities} as follows. We let
\[\lambda_c:=\inf\{\lambda\,:\,N_C>0 \mbox{ a.s.}\},\mbox{
}\lambda_u=\inf\{\lambda\,:\,N_C=1 \mbox{ a.s.}\},\]
\[\lambda_c^{*}=\sup\{\lambda\,:\,N_V>0 \mbox{ a.s. }\},\mbox{
}\lambda_u^{*}=\sup\{\lambda\,:\,N_V=1 \mbox{ a.s. }\}.\]

Our main result in $\Hyp$ is:
\begin{theorem}
\label{maintheorem} For the Poisson Boolean model with fixed radius
in $\Hyp$
\[0<\lambda_c<\lambda_u<\infty.\] Furthermore, with probability $1$,
\[(N_C,N_V)=\left\{\begin{array}{ll} (0,1), & \lambda\in[0,\lambda_c] \\ (\infty,\infty), & \lambda\in(\lambda_c,\lambda_u) \\ (1,0), & \lambda\in[\lambda_u,\infty)
  \end{array} \right.  \]
\end{theorem}
The main result in $\Hypn$ for any $n\ge 3$ is:
\begin{theorem}
\label{maintheorem2} For the Poisson Boolean model with big enough
fixed radius $R$ in $\Hypn$, $\lambda_c<\lambda_u$.
\end{theorem}

In what follows, we present several quite basic results. The proofs
of the following two lemmas, which give the possible values of $N_C$
and $N_V$ are the same as in the ${\mathbb R}^n$ case, see
Propositions 3.3 and 4.2 in \cite{meester}, and are therefore
omitted.

\begin{lemma}
\label{ergodcover} $N_C$ is an almost sure constant which equals
$0$, $1$ or $\infty$.
\end{lemma}
\begin{lemma}
\label{ergodvakant} $N_V$ is an almost sure constant which equals
$0$, $1$ or $\infty$.
\end{lemma}

Next we present some results concerning $\lambda_c$ and
$\lambda_c^{*}$.
\begin{lemma}\label{critboundlemma1}
For the Poisson Boolean model with balls of radius $R$ in $\Hypn$ it
is the case that $\lambda_c(R)>\mu(S(0,2R))^{-1}$.
\end{lemma}
The proof is identical to the ${\mathbb R}^n$ case, see Theorem 3.2 in
\cite{meester}.
\begin{prop}\label{critboundlemma2}
  Consider the Poisson Boolean model with balls of radius $R$ in
  $\Hypn$. There is $R_0<\infty$ and a constant $K=K(n)>0$ independent
  of $R$ such that for all $R\ge R_0$ we have $\lambda_c(R)\le K
  \mu(S(0,2 R))^{-1}.$
\end{prop}
\begin{proof}
  We prove the proposition using a supercritical branching process,
  the individuals of which are points in $\Hypn$. The
  construction of this branching process is done by randomly
  distorting a regular tree embedded in the space.

  Without loss of generality we assume that there is a ball centered
  at the origin, and the origin is taken to be the $0$'th generation.
  Let $a$ be such that a six-regular tree with edge length $a$ can be
  embedded in $\Hyp$ in such a way that the angles between edges at
  each vertex all equal $\pi/3$, and $d(u,v)\ge a$ for all vertices
  $u$ and $v$ in the tree. Suppose $R$ is so large that $2R-1>a$.

  Next pick three points $x_1,\, x_2,\, x_3$ on $\partial S(0,2 R)\cap \Hyp$
  such that the angles between the geodesics between the origin and
  the points is $2 \pi /3$. We define a cell associated to $x_i$ as
  the region in $S(0,2R)\setminus S(0,2R-1)$ which can be reached by a
  geodesic from the origin which diverts from the geodesic from the
  origin to $x_i$ by an angle of at most $\pi/6$.

  For every cell that contains a Poisson point, we pick one of these
  uniformly at random, and take these points to be the individuals of
  the first generation. We continue building the branching process in
  this manner. Given an individual $y$ in the n:th generation, we
  consider an arbitrary hyperbolic plane containing $y$ and its
  parent, and pick two points at distance $2R$ from $y$ in this plane
  such that the angles between the geodesics from $y$ to these two
  points and the geodesic from $y$ to its parent are all equal to
  $2\pi/3$.  Then to each of the new points, we associate a cell as
  before, and check if there are any Poisson points in them. If so, one
  is picked uniformly at random from each cell, and these points are
  the children of $y$.

  We now verify that all the cells in which the individuals of the
  branching process were found are disjoint. By construction, if $y$
  is an individual in the branching process, the angles between the
  geodesics from $y$ to its two possible children and its parent are
  all in the interval $(\pi/3,\pi)$, and therefore greater than the
  angles in a six-regular tree. Also, the lengths of these geodesics
  are in the interval $(2R-1,2R)$ and therefore larger than $a$. Thus
  by the choice of $a$, if all the individuals were in the same
  hyperbolic plane, the cells would all be disjoint.

  Suppose all individuals are in $\Hyp$, with the first individual at
  the origin. For each child of the origin we may pick two geodesics
  from the origin to infinity with angle $\theta$ less than $\pi /3$
  between them that define a sector which contains the child and all
  of its descendants and no other individuals, and the angle between
  any of these two geodesics and the geodesic between the origin and
  the child is $\theta/2$. In the same way, for each child the
  grandchildren and their corresponding descendants can be divided
  into sectors with infinite geodesics emanating from the child and so
  on. Now, such a sector emanating from an individual will contain all
  the sectors that emanates from descendants in it.

  From a sector emanating from an individual, we get a $n$-dimensional
  sector by rotating it along the geodesic going through the
  individual and its corresponding child. Then this $n$-dimensional
  sector will contain the corresponding $n$-dimensional sectors
  emanating from the child.  From this it follows that the cells will
  always be disjoint.

  Now, if the probability that a cell contains a poisson point is
  greater than $1/2$, then the expected number of children to an
  individual is greater than $1$ and so there is a positive
  probability that the branching process will never die out, which in
  turn implies that there is an unbounded connected component in the
  covered region of $\Hypn$.

  Let $B_R$ denote a cell. By \ref{iso2} there is $K_1>0$ independent
  of $R$ such that $\mu(B_R)\ge K_1 \mu(S(0,2R))$.  By the above it
  follows that \[\lambda_c(R)\le
  \frac{\log{2}}{\mu(B_R)}\le\frac{\log{2}}{K_1\mu(S(0,2R))},\]
  completing the proof.
\end{proof}
\begin{lemma}
\label{lambdavlemma2} For the Poisson Boolean model in $\Hyp$,
$\lambda_c^{*}<\infty$.
\end{lemma}
\noindent
\begin{proof}
   Let $\Gamma$ be a regular
  tiling of $\Hyp$ into congruent polygons of finite diameter. The
  polygons of $\Gamma$ can be identified with the vertices of a planar
  nonamenable transitive graph $G=(V,E)$.  Next, we define a Bernoulli
  site percolation $\omega$ on $G$. We declare each vertex $v\in V$ to
  be in $\omega$ if and only if its corresponding polygon $\Gamma(v)$
  is not completely covered by $C[\Gamma(v)]$.  Clearly, the vertices
  are declared to be in $\omega$ or not with the same probability and
  independently of each other. Now for any $v$,
  \[\lim_{\lambda\rightarrow\infty}\fatp[v \mbox{ is in $\omega$}]=0.\] Thus,
  by Theorem \ref{blps}, for $\lambda$ large enough, there are no
  infinite clusters in $\omega$. But if there are no infinite clusters
  in $\omega$, there are no unbounded components of $V$. Thus
  $\lambda_c^{*}<\infty$.
  \end{proof}

In $\Hyp$, we will need a correlation inequality for
\emph{increasing} and \emph{decreasing} events. If $\omega$ and
$\omega'$ are two realizations of a Poisson Boolean model we write
$\omega \preceq \omega'$ if any ball present in $\omega$ is also
present in $\omega'$. An event $A$ measurable with respect to the
Poisson process is said to be \emph{increasing} (respectively
\emph{decreasing}) if $\omega\preceq\omega'$ implies $1_A(\omega)\le
1_A(\omega')$ (respectively $1_A(\omega)\ge 1_A(\omega')$).
\begin{theorem}
  \label{fkg}({\sc FKG inequality}) If $A$ and $B$ are both increasing
  or both decreasing events measurable with respect to the Poisson
  process $X$, then $\fatp[A\cap B]\ge \fatp[A]\fatp[B]$.
\end{theorem}
The proof is almost identical to the proof in the ${\mathbb R}^n$
case, see Theorem 2.2 in \cite{meester}.  In particular, we will use
the following simple corollary to Theorem \ref{fkg}, the proof of
which can be found in \cite{grimmett}, which says that if
$A_1,\,A_2,\,...,\,A_m$ are increasing events with the same
probability, then
\[\fatp[A_1]\ge 1-\left(1-\fatp[\cup_{i=1}^{m}A_i]\right)^{1/m}.\] The
same holds when $A_1,\,A_2,\,...,\,A_m$ are decreasing.

For the proof of Theorem \ref{maintheorem} we need the following
lemma, the proof of which is identical to the discrete case, see
\cite{jonasson2}.

\begin{lemma}
\label{sammalemma} If
$\lim_{d(u,v)\rightarrow\infty}\fatp[u\leftrightarrow v]=0$ then
there is $a.s.$ not a unique unbounded component in $C$.
\end{lemma}

\section{The number of unbounded components in $\Hyp$}
\label{nochyp}

The aim of this section is to prove Theorem \ref{maintheorem}. We
perform the proof in the case $R=1$ but the arguments are the same for
any $R$. We first determine the possible values of $(N_C,N_V)$ for the
model in $\Hyp$. The first lemma is an application of the mass
transport principle. First, some notation is needed.
\begin{defin}
  If $H$ is a random subset of $\Hyp$ which is measurable with respect
  to the Poisson process, we say that the distribution of $H$ is
  \emph{Isom($\Hyp$)-invariant} if $g H$ has the same distribution as
  $H$ for all $g\in$\emph{ Isom($\Hyp$)}.
\end{defin}
In our applications, $H$ will typically be a union of components from
$C$ or $V$ or something similar.
\begin{lemma}\label{masslemma}
  Suppose $H$ is a random subset of $\Hyp$ which is measurable with
  respect to the Poisson process, such that its distribution is
  Isom($\Hyp$)-invariant. Also suppose that if $B$ is a bounded subset
  of $\Hyp$, then $L(B\cap\partial H)<\infty$ a.s. and $B$ intersects
  only finitely many components of $H$ a.s.  If $H$ contains only
  finite components a.s., then for any measurable $A\subset \Hyp$
  \[\E[\mu(A\cap H)]\le \E[L(A \cap
  \partial H)].\]
\end{lemma}
Before the proof we describe the intuition behind it: we place mass of
unit density in all of $\Hyp$. Then, if $h$ is a component of $H$, the
mass inside $h$ is transported to the boundary of $h$.  Then we use
the mass transport principle: the expected amount of mass transported
out of a subset $A$ equals the expected amount of mass transported
into it. Finally we combine this with the isoperimetric inequality
(\ref{isoperimetricinequality}).

\noindent
\begin{proof}
  For $A,\,B\subset\Hyp$, let \[\eta(A\times B,\,H):=\sum_h
  \frac{\mu(B\cap h)\,L(A\cap \partial h)}{L(\partial h)}.\] and let
  $\nu(A\times B):=\E[\eta(A\times B,\,H)]$. (Note that only
  components $h$ that intersect both $A$ and $B$ give a non-zero
  contribution to the sum above.) Since the distribution of $H$ is
  Isom($\Hyp$)-invariant, we get for each $g\in $Isom($\Hyp$)
\[\nu(gA\times g B)=\E[\eta(g A\times g B,\,H)]=\E[\eta(g A\times g B,\,g
H)]\]\[=\E[\eta(A\times B,\,H)]=\nu(A\times B).\] Thus, $\nu$ is a
diagonally invariant positive measure on $\Hyp\times \Hyp$.  We have
$\nu(\Hyp\times A)=\E\left[\mu(A\cap H)\right]$ and \[\nu(A\times
\Hyp)=\E\left[\sum_h
  \frac{\mu(h)\,L(A\cap\partial h)}{L(\partial h)}\right]\le
\E[L(A\cap\partial H)]\] where the last inequality follows from the
linear isoperimetric inequality. Hence, the claim follows by Theorem
\ref{masstransport}.
\end{proof}

In the following lemmas, we exclude certain combinations of $N_C$ and
$N_V$. The first lemma can be considered as a continuous analogue to
Lemma 3.3 in \cite{itai1}.
\begin{lemma}
\label{lemma00}
If $H$ is a union of components from ${\mathcal C}$ and ${\mathcal V}$
such that the distribution of $H$ is Isom($\Hyp$)-invariant, then $H$
and/or $H^c$ contains unbounded components almost surely.
\end{lemma}
\begin{proof}
  Suppose $H$ and $D:=H^c$ contains only finite components, and let in
  this proof ${\mathcal H}_0$ and ${\mathcal D}_0$ be the collections
  of the components of $H$ and $D$ respectively.  Then every element
  $h$ of ${\mathcal H}_0$ is surrounded by a unique element $h'$ of
  ${\mathcal D}_{0}$, which in turn is surrounded by a unique element
  $h''$ of ${\mathcal H}_0$. In the same way, every element $d$ of
  ${\mathcal D}_{0}$ is surrounded by a unique element $d'$ of
  ${\mathcal H}_0$ which in turn is surrounded by a unique element
  $d''$ of ${\mathcal D}_0$. Inductively, for $j\in {\mathbb N}$, let
  ${\mathcal H}_{j+1}:=\{h''\,:\,h\in{\mathcal H}_{j}\}$ and
  ${\mathcal D}_{j+1}:=\{d''\,:\,d\in {\mathcal D}_{j}\}$. Next, for
  $r\in {\mathbb N}$, let
  \[A_{r}:=\bigcup_{j=0}^{r}(\{h\in {\mathcal
    H}_0\,:\,\sup\{i\,:\,h\in{\mathcal H}_{i}\}=j\}\cup\{d\in
  {\mathcal D}_0\,:\,\sup\{i\,:\,d\in{\mathcal D}_{i}\}=j\}).\] In
  words, ${\mathcal H}_j$ and ${\mathcal D}_j$ define layers of
  components from $H$ and $D$. Thus $A_r$ is the union of all layers
  of components from $H$ and $D$ that have at most $r$ layers inside
  of them.  Now let $B$ be some ball in $\Hyp$. Note that $L(B\cap
  \partial A_r)\le L(B\cap \partial C)$ and $\E[L(B\cap \partial
  C)]<\infty$.  Also, almost surely, there is some random $r_0$ such
  that $B$ will be completely covered by $A_r$ for all $r\ge r_0$.
  Thus the dominated convergence theorem gives
  \[\lim_{r\rightarrow \infty}\E[\mu(B\cap A_r)]=\mu(B)\mbox{ and
  }\lim_{r\rightarrow \infty}\E[L(B \cap \partial A_r)]=0.\] Since the
  distribution of $A_r$ is Isom($\Hyp$)-invariant we get by Lemma
  \ref{masslemma} that there is $r_1<\infty$ such that for $r\ge r_1$,
  \[{\bf P}[A_r \mbox{ has unbounded components}]>0.\] But by
  construction, for any $r$ it is the case that $A_r$ has only finite
  components. Hence the initial assumption is false.\end{proof}
\begin{lemma}
\label{lemma1infty}
The cases $(N_C,N_V)=(\infty,1)$ and $(N_C,N_V)=(1,\infty)$ have
probability $0$.
\end{lemma}
\noindent
\begin{proof} Suppose $N_C=\infty$. First we show that it is possible
  to pick $r>0$ such that the event \begin{multline*}A(x,r):=
    \\\{S(x,r) \mbox{ intersects at least } 2 \mbox{ disjoint
      unbounded components of } C[S(x,r)^c]\}\end{multline*} has
  positive probability for $x\in \Hyp$. Suppose $S(x,r_0)$ intersects
  an unbounded component of $C$ for some $r_0>0$.  Then if $S(x,r_0)$
  does not intersect some unbounded component of $C[S(x,r_0)^c]$,
  there must be some ball centered in $S(x,r_0+2)\backslash
  S(x,r_0+1)$ being part of an unbounded component of
  $C[S(x,r_0+1)^c]$, which is to say that $S(x,r_0+1)$ intersects an
  unbounded component of $C[S(x,r_0+1)^c]$.  Clearly we can find
  $\tilde{r}$ such that
  \[B(x,\tilde{r}):=\{S(x,\tilde{r}) \mbox{ intersects at least }
  3 \mbox{ disjoint unbounded components of } C\}.\] By the above
  discussion it follows that $\fatp[A(x,\tilde{r})\cup
  A(x,\tilde{r}+1)]>0$, which proves the existence of $r$ such that
  $A(x,r)$ has positive probability. Pick such an $r$ and let
  $E(x,r):=\{S(x,r)\subset C[S(x,r)]\}$. $E$ has positive probability
  and is independent of $A$ so $A\cap E$ has positive probability. By
  planarity, on $A\cap E$, $V$ contains at least 2 unbounded
  components.  So with positive probability, $N_V>1$. By Lemma
  \ref{ergodvakant}, $N_V=\infty$ a.s. This finishes the first part of
  the proof.
Now instead suppose $N_V=\infty $ and pick $r>0$ such that
\[A(x,r):=\{S(x,r) \mbox{ intersects at least two unbounded components
  of } V\}\] has positive probability. Let
\[B(x,r):=\{C[S(x,r+1)^c] \mbox{ contains at least 2 unbounded
  components}\}.\] On $A$, $C\backslash S(x,r)$ contains at least two
unbounded components, which in turn implies that $B$ occurs.  Since
$\fatp[A]>0$ this gives $\fatp[B]>0$. Since $B$ is independent of
$F(x,r):=\{|X(S(x,r+1))|=0\}$ which has positive probability,
$\fatp[B\cap F]>0$. On $B\cap F$, $C$ contains at least two unbounded
components. By Lemma \ref{ergodcover} we get $N_C=\infty$ a.s.
\end{proof}

The proof of the next lemma is very similar to the discrete case, see
Lemma 11.12 in \cite{grimmett}, but is included for the convenience of
the reader.
\begin{lemma}
\label{lemma11}
The case $(N_C,N_V)=(1,1)$ has probability $0$.
\end{lemma}
\begin{proof}
  Assume $(N_C,N_V)=(1,1)$ a.s.  Fix $x\in \Hyp$. Denote by
  $A^{u}_{C}(k)$ (respectively $A^d_{C}(k)$, $A^r_{C}(k)$,
  $A^l_{C}(k)$) the event that the uppermost (respectively lowermost,
  rightmost, leftmost) quarter of $\partial S(x,k)$ intersects an
  unbounded component of $C\backslash S(x,k)$. Clearly, these events
  are increasing. Since $N_C=1$ a.s.,
\[
\lim_{k\rightarrow \infty}\fatp[A^{u}_C(k)\cup A^d_C(k)\cup
A^r_C(k)\cup A^l_C(k)]=1.
\]
Hence by the corollary to the FKG-inequality, $\lim_{k\rightarrow
  \infty}\fatp[A^t_C(k)]=1$ for \\$t\in\{u,\,d,\,r,\,l\}$. Now let
$A^u_V(k)$ (respectively $A^l_V(k)$, $A^r_V(k)$, $A^l_V(k)$) be the
event that the uppermost (respectively lowermost, rightmost, leftmost)
quarter of $\partial S(x,k)$ intersects an unbounded component of
$V\backslash S(x,k)$. Since these events are decreasing, we get in the
same way as above that $\lim_{k\rightarrow\infty}\fatp[A_V^t(k)]=1$
for $t\in\{u,\,d,\,r,\,l\}$. Thus we may pick $k_0$ so big that
$\fatp[A_C^t(k_0)]>7/8$ and $\fatp[A_V^t(k_0)]>7/8$ for
$t\in\{u,\,d,\,r,\,l\}$. Let
\[A:=A_C^u(k_0)\cap A_C^d(k_0)\cap A_V^l(k_0)\cap A_V^r(k_0).\]
Bonferroni's inequality implies $\fatp[A]>1/2$. On $A$, $C\backslash
S(x,k_0)$ contains two disjoint unbounded components. Since $N_C=1$
a.s., these two components must almost surely on $A$ be connected. The
existence of such a connection implies that there are at least two
unbounded components of $V$, an event with probability $0$. This gives
$\fatp[A]=0$, a contradiction.
\end{proof}
\begin{prop}
\label{antalsproposition}
Almost surely, $(N_C,N_V)\in\{(1,0),\,(0,1),\,(\infty,\infty)\}$.
\end{prop}
\noindent
\begin{proof}
  By Lemmas \ref{ergodcover} and \ref{ergodvakant}, each of $N_C$ and
  $N_V$ is in $\{0,\,1,\,\infty\}$. Lemma \ref{lemma00} with $H\equiv
  C$ rules out the case $(0,0)$. Hence Lemmas \ref{lemma1infty} and
  \ref{lemma11} imply that it remains only to rule out the cases
  $(0,\infty)$ and $(\infty,0)$. But since every two unbounded
  components of $C$ must be separated by some unbounded component of
  $V$, $(\infty,0)$ is impossible. In the same way, $(0,\infty)$ is
  impossible.
\end{proof}

\subsection{The situation at $\lambda_c$ and $\lambda_c^*$}

It turns out that to prove the main theorem, it is necessary to
investigate what happens regarding $N_C$ and $N_V$ at the intensities
$\lambda_c$ and $\lambda_{c}^{*}$. Our proofs are inspired by the
proof of Theorem 1.1 in \cite{itai2}, which says that critical
Bernoulli bond and site percolation on nonamenable Cayley graphs does
not contain infinite clusters.

\begin{theorem}
\label{lambdactheorem}
At $\lambda_c$, $N_C=0$ a.s.
\end{theorem}

\noindent
\begin{proof}
  We begin with ruling out the possibility of a unique unbounded
  component of $C$ at $\lambda_c$.  Suppose $\lambda=\lambda_c$ and
  that $N_C=1$ a.s. Denote the unique unbounded component of $C$ by
  $U$. By Proposition \ref{antalsproposition}, $V$ contains only
  finite components a.s. Let $\epsilon>0$ be small and remove each
  point in $X$ with probability $\epsilon$ and denote by
  $X_{\epsilon}$ the remaining points. Furthermore, let
  $C_{\epsilon}=\cup_{x\in X_{\epsilon}}S(x,1)$.  Since $X_{\epsilon}$
  is a Poisson process with intensity $\lambda_c(1-\epsilon)$ it
  follows that $C_{\epsilon}$ will contain only bounded components
  a.s. Let ${\mathcal C}_{\epsilon}$ be the collection of all
  components of $C_{\epsilon}$. We will now construct $H_{\epsilon}$
  as a union of elements from ${\mathcal C}_{\epsilon}$ and ${\mathcal
    V}$ such that the distribution of $H_{\epsilon}$ will be
  Isom($\Hyp$)-invariant.  For each $z\in\Hyp$ we let
  $U_{\epsilon}(z)$ be the union of the components of $U\cap
  C_{\epsilon}$ being closest to $z$. We let each $h$ from ${\mathcal
    C}_{\epsilon}\cup{\mathcal V}$ be in $H_{\epsilon}$ if and only if
  $\sup_{z\in h}d(z,U)<1/\epsilon$ and
  $U_{\epsilon}(x)=U_{\epsilon}(y)$ for all $x,\,y\in h$.
  We want to show that for $\epsilon$ small enough,
  $H_{\epsilon}$ contains unbounded components with positive
  probability. Let $B$ be some ball. It is clear that
  $L(B \cap \partial H_{\epsilon})\rightarrow 0$ a.s. and also that $\mu(B\cap H_{\epsilon})\rightarrow \mu(B)$ a.s. when $\epsilon\rightarrow
  0$. Also $L(B\cap \partial H_{\epsilon}) \le L(B \cap (\partial
  C_{\epsilon}\cup \partial C))$ and $\E[L(B \cap (\partial
  C_{\epsilon}\cup \partial C))]\le K<\infty$ for some constant $K$
  independent of $\epsilon$. By the dominated convergence theorem,
  we have
  \[\lim_{\epsilon\rightarrow 0}\E[\mu(B\cap
  H_{\epsilon})]=\mu(B)\mbox{ and }\lim_{\epsilon\rightarrow
    0}\E[L(B\cap\partial H_{\epsilon})]=0.\]
\noindent
Therefore we get by Lemma \ref{masslemma} that $H_{\epsilon}$ contains
unbounded components with positive probability when $\epsilon$ is
small enough. Suppose $h_1,h_2,...$ is an infinite sequence of
distinct elements from ${\mathcal C}_{\epsilon}\cup{\mathcal V}$ such
that they constitute an unbounded component of $H_{\epsilon}$.  Then
$U_{\epsilon}(x)=U_{\epsilon}(y)$ for all $x,y$ in this component.
Hence $U\cap C_{\epsilon}$ contains an unbounded component (this
particular conclusion could not have been made without the condition
$\sup_{z\in h}d(z,U)<1/\epsilon$ in the definition of
$U_{\epsilon}(z)$).  Therefore we conclude that the existence of an
unbounded component in $H_{\epsilon}$ implies the existence of an
unbounded component in $C_{\epsilon}$. Hence $C_{\epsilon}$ contains
an unbounded component with positive probability, a contradiction.

We move on to rule out the case of infinitely many unbounded
components of $C$ at $\lambda_c$. Assume $N_C=\infty$ a.s. at
$\lambda_c$. As in the proof of Lemma \ref{lemma1infty}, we choose $r$
such that for $x\in \Hyp$ the event
\begin{multline*}A(x,r):=\\\{S(x,r) \mbox{ intersects at least } 3 \mbox{ disjoint unbounded
  components of } C[S(x,r)^c]\}\end{multline*} has positive probability.  Let
$B(x,r):=\{S(x,r)\subset C[S(x,r)]\}$ for $x\in \Hyp$. Since $A$ and
$B$ are independent, it follows that $A\cap B$ has positive
probability. On $A\cap B$, $x$ is contained in an unbounded component
$U$ of $C$. Furthermore, $U\backslash S(x,r+1)$ contains at least
three disjoint unbounded components. Now let $Y$ be a Poisson process
independent of $X$ with some positive intensity. We call a point $y\in
\Hyp$ a \emph{encounter point} if
\begin{itemize}
\item $y\in Y$;
\item $A(y,r)\cap B(y,r)$ occurs;
\item $S(y,2(r+1))\cap Y=\{y\}$.
\end{itemize} The third condition above means that if $y_1$ and $y_2$
are two encounter points, then $S(y_1,r+1)$ and $S(y_2,r+1)$ are
disjoint sets. By the above, it is clear that given $y\in Y$, the
probability that $y$ is an encounter point is positive.  

We now move on to show that if $y$ is an encounter point and $U$ is
the unbounded component of $C$ containing $y$, then each of the
disjoint unbounded components of $U\backslash S(y,r+1)$ contains a
further encounter point.

Let $m(s,t)=1$ if $t$ is the unique encounter point closest to $s$ in
$C$, and $m(s,t)=0$ otherwise. Then let for measurable sets
$A,\,B\subset\Hyp$
\[\eta(A\times B)=\sum_{s\in Y(A)}\sum_{t\in Y(B)}m(s,t)\]
and
\[\nu(A\times B)=\E[\eta(A\times B)].\]
Clearly, $\nu$ is a positive diagonally invariant measure on
$\Hyp\times \Hyp$.  Suppose $A$ is some ball in $\Hyp$. Since
$\sum_{t\in Y}m(s,t)\le 1$ we get $\nu(A\times \Hyp)\le
\E[|Y(A)|]<\infty$. On the other hand, if $y$ is an encounter point
lying in $A$ and with positive probability there is no encounter point
in some of the unbounded components of $U\backslash S(y,r+1)$ we get
$\sum_{s\in Y}\sum_{t\in Y(A)}m(s,t)=\infty$ with positive
probability, so $\nu(\Hyp\times A)=\infty$, which contradicts Theorem
\ref{masstransport}.

The proof now continues with the construction of a forest $F$, that is
a graph without loops or cycles. Denote the set of encounter points by
$T$, which is a.s. infinite by the above. We let each $t\in T$ represent a
vertex $v(t)$ in $F$. For a given $t\in T$, let $U(t)$ be the unbounded
component of $C$ containing $t$. Then let $k$ be the number of
unbounded components of $U(t)\backslash S(t,r+1)$ and denote these
unbounded components by $C_1$, $C_2$,..., $C_k$. For
$i=1,\,2,\,...,\,k$ put an edge between $v(t)$ and the vertex
corresponding to the encounter point in $C_i$ which is closest to $t$
in $C$ (this encounter point is unique by the nature of the Poisson
process).

Next, we verify that $F$ constructed as above is indeed a forest. If
$v$ is a vertex in $F$, denote by $t(v)$ the encounter point
corresponding to it. Suppose $v_0,\,v_1,...,\,v_n=v_0$ is a cycle of
length $\ge 3$, and that $d_C(t(v_0),t(v_1))<d_C(t(v_1),t(v_2))$. Then
by the construction of $F$ it follows that
$d_C(t(v_1),t(v_2))<d_C(t(v_2),t(v_3))<...<d_C(t(v_{n-1}),t(v_0))<d_C(t(v_0),t(v_1))$
which is impossible. Thus we must have that $d_C(t(v_i),t(v_{i+1}))$
is the same for all $i\in\{0,1,..,n-1\}$. The assumption
$d_C(t(v_0),t(v_1))>d_C(t(v_1),t(v_2))$ obviously leads to the same
conclusion.  But if $y\in Y$, the probability that there are two other
points in $Y$ on the same distance in $C$ to $y$ is $0$.  Hence,
cycles exist with probability $0$, and therefore $F$ is almost surely
a forest.

Now define a bond percolation $F_{\epsilon}\subset F:$ Define
$C_{\epsilon}$ in the same way as above.  Let each edge in $F$ be in
$F_{\epsilon}$ if and only if both encounter points corresponding to
its end-vertices are in the same component of $C_{\epsilon}$. Since
$C_{\epsilon}$ contains only bounded components, $F_{\epsilon}$
contains only finite connected components.

For any vertex $v$ in $F$ we let $K(v)$ denote the connected component
of $v$ in $F_{\epsilon}$ and let $\partial_{F} K(v)$ denote the inner
vertex boundary of $K(v)$ in $F$. Since the degree of each vertex in
$F$ is at least $3$, and $F$ is a forest, it follows that at least
half of the vertices in $K(v)$ are also in $\partial_F K(v)$. Thus we
conclude
\[\fatp[x\in T,\, v(x)\in \partial_F K(v(x))|x\in Y]\ge \frac{1}{2}\fatp[x\in T|x\in Y].\]
The right-hand side of the above is positive and independent of
$\epsilon$.  But the left-hand side tends to $0$ as $\epsilon$ tends
to $0$, since when $\epsilon$ is small, it is unlikely that an edge in
$F$ is not in $F_{\epsilon}$. This is a contradiction.\end{proof}

By Proposition \ref{antalsproposition}, if $N_C=0$ a.s., then
$N_V=1$ a.s. Thus we have an immediate corollary to Theorem
\ref{lambdactheorem}.

\begin{corollary}
At $\lambda_c$, $N_V=1$ a.s.
\end{corollary}

Next, we show the corresponding results for $\lambda_c^{*}$.
Obviously, the nature of $V$ is quite different from that of $C$, but
still the proof of Theorem \ref{lambdavtheorem} below differs only in
details to that of Theorem \ref{lambdactheorem}. We include it for the
convenience of the reader.

\begin{theorem}
\label{lambdavtheorem} At $\lambda_c^{*}$, $N_V=0$ a.s.
\end{theorem}
\noindent
\begin{proof}
  Suppose $N_V=1$ a.s. at $\lambda_c^{*}$ and denote the unbounded
  component of $V$ by $U$. Then $C$ contains only finite components a.s. by Proposition
  \ref{antalsproposition}. Let $\epsilon>0$ and let $Z$ be a Poisson
  process independent of $X$ with intensity $\epsilon$.  Let
  $C_{\epsilon}:=\cup_{x\in X\cup Z} S(x,1)$ and
  $V_{\epsilon}:=C_{\epsilon}^c$. Since $X\cup Z$ is a Poisson process
  with intensity $\lambda_c^{*}+\epsilon$ it follows that $C_{\epsilon}$
  has a unique unbounded component a.s. and hence $V_{\epsilon}$
  contains only bounded components a.s. Let ${\mathcal V}_{\epsilon}$
  be the collection of all components of $V_{\epsilon}$. Define
  $H_{\epsilon}$ in the following way: For each $z\in\Hyp$ we let
  $U_{\epsilon}(z)$ be the union of the components of $U\cap
  V_{\epsilon}$ being closest to $z$. We let each $h\in{\mathcal
    C}\cup{\mathcal V}_{\epsilon}$ be in $H_{\epsilon}$ if and only if
  $\sup_{z\in h}d(z,U)<1/\epsilon$ and
  $U_{\epsilon}(x)=U_{\epsilon}(y)$ for all $x,\,y\in h$.  As in the
  proof of Theorem \ref{lambdactheorem}, for $\epsilon>0$ small
  enough, $H_{\epsilon}$ contains an unbounded component with positive
  probability, and therefore $V_{\epsilon}$ contains an unbounded
  component with positive probability, a contradiction.

  Now suppose that $N_V=\infty$ a.s. at $\lambda_c^{*}$. Then also
  $N_C=\infty$ by Proposition \ref{antalsproposition}. Therefore, for
  $x\in \Hyp$, we can choose $r>1$ big such that the intersection of
  the two independent events \begin{multline*}A(x,r):=\\\{S(x,r)
    \mbox{ intersects at least } 3 \mbox{ disjoint unbounded
      components of } C[S(x,r)^c]\}\end{multline*} and
  $B(x,r):=\{|X(S(x,r))|=0\}$ has positive probability. Next, suppose
  that $Y$ is a Poisson process independent of $X$ with some positive
  intensity. Now we redefine what an encounter point is: call $y\in
  \Hyp$ an encounter point if
\begin{itemize}
\item $y\in Y$;
\item $A(y,r)\cap B(y,r)$ occurs;
\item $S(y,2r)\cap Y=\{y\}$.
\end{itemize}
By the above discussion, \[\fatp[y \mbox{ is an encounter point
}|\,y\in Y]>0.\] If $y$ is a encounter point, $y$ is contained in an
unbounded component $U$ of $V$ and $U\backslash S(y,r)$ contains at
least $3$ disjoint unbounded components.  Again we construct a
forest $F$ using the encounter points and define a bond percolation
$F_{\epsilon}\subset F$. Let $V_{\epsilon}$ be defined as above.
Each edge of $F$ is declared to be in $F_{\epsilon}$ if and only if
both its end-vertices are in the same component of $V_{\epsilon}$.
The proof is now finished in the same way as Theorem
\ref{lambdactheorem}.\end{proof}

Again, Proposition \ref{antalsproposition} immediately implies the
following corollary:
\begin{corollary}
  At $\lambda_c^{*}$, $N_C=1$ a.s.
\end{corollary}
\subsection{Proof of Theorem \ref{maintheorem} }
Here we combine the results from the previous sections to prove our
main theorem in $\Hyp$.

{\em Proof of Theorem \ref{maintheorem}.} If
$\lambda<\lambda_u$ then Proposition \ref{antalsproposition} implies
$N_V>0$ a.s. giving $\lambda\le \lambda_{c}^{*}$. If
$\lambda>\lambda_{u}$ the same proposition gives $N_V=0$ a.s. giving
$\lambda\ge \lambda_{c}^{*}$.  Thus
\begin{equation}
\label{finaleq1}
\lambda_{u}=\lambda_{c}^{*}.
\end{equation}
By Theorem \ref{lambdactheorem} $N_C=0$ a.s. at $\lambda_c$, so
$N_V>0$ a.s. at $\lambda_{c}$ by Proposition \ref{antalsproposition}.
Thus by Theorem \ref{lambdavtheorem}
\begin{equation}
\label{finaleq2}
\lambda_c<\lambda_{c}^{*}.
\end{equation}
Hence the desired conclusion follows by (\ref{finaleq1}),
(\ref{finaleq2}) and Lemma \ref{lambdavlemma2}. \qed

\section{The number of unbounded components in  $\Hypn$}
\label{nochypn}

This section is devoted to the proof of Theorem \ref{maintheorem2}.

 \noindent{\em First part of proof of Theorem \ref{maintheorem2}.}
In view of Lemma \ref{sammalemma}, it is enough to show that
$\fatp[u\leftrightarrow v]\rightarrow 0$ as $d(u,v)\rightarrow
\infty$ for some intensity above $\lambda_c$.  We use a duplication
trick.  Let $X_1$ and $X_2$ be two independent copies of the Poisson
Boolean model. If we for some $\epsilon>0$ can find points $u$ and
$v$ on an arbitrarily large distance from each other such that $u$
is connected to $v$ in $X_1$ with probability at least $\epsilon$,
then the event
\[B(u,v):=\{u \mbox{ is connected to } v \mbox{ in both }X_1\mbox{ and
}X_2\}\] has probability at least $\epsilon^2$. So it is enough to
show that $\fatp[B(u,v)]\rightarrow 0$ as $d(u,v)\rightarrow \infty$
at some intensity above $\lambda_c$.

Fix points $u$ and $v$ and suppose $d(u,v)=d$. Let $k=\lceil d/(2
R)\rceil$. That is, $k$ is the smallest number of balls of radius
$R$ needed to connect the points $u$ and $v$. Thus, for $B(u,v)$ to
occur, there must be at least one sequence of at least $k$ distinct
connected balls in $X_1$, such that the first ball contains $u$ and
the last ball contains $v$, and at least one such sequence of balls
in $X_2$. This in turn obviously implies that there is at least one
sequence of at least $k$ connected balls in $X_1$ such that the
first ball contains $u$, and the last ball intersects the first ball
of a sequence of at least $k$ connected balls in $X_2$, where the
last ball in this sequence contains $u$. In this sequence of at
least $2 k$ balls, the center of the first ball is at distance at
most $2R$ from the center of the last ball.

Let $l\ge 2 k$. Next we estimate the expected number of sequences of
balls as above of length $l$. Denote this number by $N(l)$. Now, if
we consider sequences of balls as above of length $l$, without the
condition that the last ball contains $u$, then the expected number
of such sequences is easily seen to be bounded by $\lambda^l \mu
(S(0,2R))^l$ (as for example in the proof of Theorem 3.2 in
\cite{meester}). Let $P_R(l)$ be the probability that the center of
the last ball in such a sequence is at most at distance $2R$ from
the center of the first ball. Then $N(l)\le \lambda^l \mu(S(0,2R))^l
P_R(l)$.

Now \[\fatp[B(u,v)]\le\sum_{l=2 k}^{\infty}N(l)\le\sum_{l=2
  k}^{\infty}\left(\lambda \mu (S(0,2R))\right)^l P_R(l).\] We will
now estimate the terms in the sum above.
\begin{lemma}
\label{brownlemma} Suppose $X_0,X_1,...X_k$ is a sequence of
distinct points in a Poisson point process in $\Hypn$ such that
$d(X_i,X_{i+1})<2 R$ for $i=0,1,...,k-1$. Then there is a sequence
of i.i.d. random variables $Y_1,\,Y_2,...$ with positive mean such
that
\[\fatp[d(X_0,X_k)\le 2R]\le \fatp[\sum_{i=1}^{k-1}Y_i\le 2R].\] In other words,
$P_R(k)\le \fatp[\sum_{i=1}^{k-1}Y_i\le 2R]$.
\end{lemma}

The distribution of $Y_i$ will be defined in the proof.\bigskip

\noindent {\em First part of proof of Lemma \ref{brownlemma}.}
  Note that given the point $X_i$, the
  distribution of the point $X_{i+1}$ is the uniform distribution on
  $S(X_i,2 R)$. Put $d_i:=d(X_i,X_{i+1})$. Then $d_0,d_1...$ is a
  sequence of independent random variables with density

\begin{equation}
  \frac{d}{dr} \frac{\mu(S(0,r))}{\mu(S(0,2R))}=\frac{\sinh(r)^{n-1}}{\int_{0}^{2R}\sinh(t)^{n-1}\,dt} \mbox{ for }r\in [0,2R].
\end{equation}

\noindent
Next we write

\begin{equation}
\label{telesumma}
\fatp[d(X_0,X_k)<2 R]=\fatp\left[\sum_{i=0}^{k-1}(d(X_0,X_{i+1})-d(X_0,X_i))<2 R\right].
\end{equation}

\noindent The terms in the sum \ref{telesumma} are neither
independent nor identically distributed. However, we will see that
the sum is always larger than a sum of i.i.d. random random
variables with positive mean.  Suppose without loss of generality
that $X_0$ is at the origin. Let $\gamma_i$ be the geodesic between
$0$ and $X_i$ and let $\varphi_i$ be the geodesic between $X_i$ and
$X_{i+1}$ for $i\ge 1$. Let $\theta_i$ be the angle between
$\gamma_i$ and $\varphi_i$ for $i\ge 1$ and let $\theta_0=\pi$. Then
$\theta_1,\theta_2,...$ is a sequence of independent random
variables, uniformly distributed on $[0,\pi]$. Since the geodesics
$\gamma_i$, $\gamma_{i+1}$ and $\varphi_i$ lie in the same
hyperbolic plane, we can express $d(0,X_{i+1})$ in terms of
$d(0,X_i)$, $d(X_i,X_{i+1})$ and $\theta_i$ using the first law of
cosines for triangles in hyperbolic space (see \cite{ratcliffe},
Theorem 3.5.3), which gives that

\begin{equation}\begin{split}
    d(0,X_{i+1})-d(0,X_i)&= \cosh^{-1}\bigg(\cosh(d_i)\cosh(d(0,X_i)) \\
    &-\sinh(d_i)\sinh(d(0,X_i))\cos(\theta_i)\bigg)-d(0,X_i).\end{split}
\end{equation}

\noindent
Next we prove a lemma that states that the random variable above
dominates a random variable which is independent of $d(0,X_i)$.  Put
\[f(x,y,\theta):=\cosh^{-1}(\cosh(x)\cosh(y)-\sinh(x)\sinh(y)\cos(\theta))-y.\]

\begin{lemma}
\label{funktionslemma}
For fixed $x$ and $\theta$, the function $f(x,y,\theta)$ is strictly
decreasing in $y$ and
$g(x,\theta):=\lim_{y\rightarrow\infty}f(x,y,\theta)=\log(\cosh(x)-\sinh(x)\cos(\theta))$.
\end{lemma}
\begin{proof}
  For simplicity write $a=a(x):=\cosh(x)$ and
  $b=b(x,\theta):=\sinh(x)\cos(\theta)$. Then by rewriting
\begin{equation}
  f(x,y,\theta)=\log\left(\frac{a\cosh(y)-b\sinh(y)+\sqrt{(a\cosh(y)-b\sinh(y))^2-1}}{\exp(y)}\right)
\end{equation}
\noindent
we get by easy calculations that the limit as $y\rightarrow\infty$ is
as desired.
\noindent
It remains to show that $f_y^{'}(x,y,\theta)<0$ for all $x,\,y$ and
$\theta$. We have that
\begin{equation}
  f_y^{'}(x,y,\theta)=-1+\frac{-b\cosh(y)+a\sinh(y)}{\sqrt{-1+a\cosh(y)-b\sinh(y)}\sqrt{1+a\cosh(y)-b\sinh(y)}}
\end{equation}
\noindent
which is less than $0$ if
\begin{equation}
\label{olikhet}
\sqrt{-1+a\cosh(y)-b\sinh(y)}\sqrt{1+a\cosh(y)-b\sinh(y)}>a\sinh(y)-b\cosh(y)
\end{equation}
\noindent
If the right hand side in \ref{olikhet} is negative then we are done,
otherwise, taking squares and simplifying gives that the inequality
\ref{olikhet} is equivalent to the simpler inequality
\[a^2-b^2>1\]
\noindent
which holds since
$a^2-b^2=\cosh^2(x)-\sinh^2(x)\cos^2(\theta)>\cosh^2(x)-\sinh^2(x)=1$,
completing the proof of the lemma.\end{proof}
\noindent {\em Second part of proof of Lemma \ref{brownlemma}. }
Letting $Y_i:=g(d_i,\theta_i)$ we have (since $Y_0>0$),
\begin{equation}
  \label{summaskattning}
  \fatp[d(X_0,X_k)<2 R]\le \fatp\left[\sum_{i=0}^{k-1}Y_i<2 R
  \right]\le \fatp\left[\sum_{i=1}^{k-1}Y_i<2 R\right]
\end{equation} where $g$ is as in Lemma \ref{funktionslemma}, which concludes the
proof.\qed\bigskip

We now want to bound the probability in Lemma
\ref{brownlemma}, and for this we have the following technical lemma,
which in a slightly different form than below is due to Patrik Albin.
\begin{lemma}\label{palbinlemma}
  Let $Y_i$ be defined as above. There is a function $h(R,\epsilon)$
  such that for any $\epsilon\in (0,1)$ we have $h(R,\epsilon)\sim A
 e^{-R(1-\epsilon)}$ as $R\rightarrow\infty$ for some constant
  $A=A(\epsilon)\in (0,\infty)$ independent of $R$ and such that for any $R>0$,
\begin{equation}
  \fatp\left[\sum_{i=1}^k Y_i<2R\right]\le h(R,\epsilon)^k e^R.
\end{equation}
\end{lemma}
\begin{proof}
  Let $K$ be the complete elliptic integral of the first kind (see
  \cite{ErdMagObeTri}, pp. 313-314). Then we have
  \[\begin{split}
    \E[e^{-Y_1/2}|d_1]&=\E\left[\frac{1}{\sqrt{\cosh(d_1)-\sinh(d_1)\cos(\theta_1)}}\bigg{|}d_1\right]\\
    &=\E\left[\frac{e^{-d_1/2}}{\sqrt{1-\cos(\theta_1/2)^2(1-e^{-2d_1})}}\bigg{|}d_1\right]\\
    &=\frac{2 e^{-d_1/2}K(\sqrt{1-e^{-2d_1}})}{\pi}.\end{split}\]
\noindent Using the relation $K(x)=\pi$ $_2F_1(1/2,1/2,1,x)/2$ where
$_2F_1$ is the hypergeometric function (see \cite{ErdMagObeTri},
Equation 13.8.5), we have
\[\E[e^{-Y_1/2}|d_1]=e^{-d_1/2}\,
_2F_1(1/2,1/2,1,1-e^{-2d_1}).\] Since $_2F_1(1/2,1/2,1,\cdot)$ is
continuous on $\{z\in {\mathbb C}\,:\,|z|\le\rho\}$ for any
$\rho\in(0,1)$, this gives
\begin{equation}\label{skattning1}
  \E[e^{-Y_1/2}|d_1]\le A_1 e^{-d_1/2}\mbox{ for }d_1\le
  x_0,\end{equation}
for some constant $A_1(x_0)>0$, for any $x_0>0$. Large values of
$d_1$ makes the argument of $_2F_1(1/2,1/2,1,1-e^{-2d_1})$
approach the radius of convergence $1$ of $_2F_1(1/2,1/2,1,\cdot)$
so we perform the quadratic transformation
\[_2F_1(a,b,2b,x)=(1-x)^{-a/2}{}_2F_1\left(a,2b-a,b+1/2,-\frac{(1-\sqrt{1-x})^2}{4\sqrt{1-x}}\right),\]
(see \cite{ErdMagObeTri0}, Equation 2.11.30), giving
\[\E[e^{-Y_1/2}|d_1]={}
_2F_1\left(1/2,1/2,1,-e^{d_1}(1-e^{-d_1})^2/4\right).\]
\noindent By the asymptotic behaviour of the hypergeometric function
(here the analytic continuation of the hypergeometric function is
used), we have
\[|_2F_1(1/2,1/2,1,x)|\sim A_2\frac{\log{|x|}}{\sqrt{|x|}}\] as
$|x|\rightarrow\infty$ (see \cite{ErdMagObeTri0}, Equation 2.3.2.9),
for some constant $A_2>0$. Combining this with \ref{skattning1} we
get
\[\E[e^{-Y_1/2}|d_1]\le A_3(1+d_1)e^{-d_1/2}\le
A_4 e^{-(1-\epsilon)d_1/2}\] \noindent for $d_1>0$, for any
$\epsilon\in(0,1)$, for some constants $A_3>0$ and $A_4(\epsilon)>0$.
Thus \[\begin{split}
  \E[e^{-Y_1/2}]&\le\E[A_4 e^{-d_1(1-\epsilon)/2}]\\
  &=A_4\frac{\int_{0}^{2R}\sinh(t)^{n-1}e^{-t(1-\epsilon)/2}\,dt}{\int_{0}^{2R}\sinh(t)^{n-1}
    \,dt}\\ \end{split}\]
Clearly $h(R,\epsilon):=A_4
\int_{0}^{2R}\sinh(t)^{n-1}e^{-t(1-\epsilon)/2}\,dt \diagup
\int_{0}^{2R}\sinh(t)^{n-1} \,dt\sim A e^{-R(1-\epsilon)}$ as
$R\rightarrow \infty$ for some constant $A\in(0,\infty)$.  Finally we
get using Markov's inequality that
\[\begin{split}\fatp\left[\sum_{i=1}^k Y_i<
    2R\right]&=\fatp\left[e^{-\frac{1}{2}\sum_{i=1}^k
      Y_i}>e^{-R}\right]\\ & \le
  e^R\E\left[e^{-\frac{1}{2}\sum_{i=1}^k
      Y_i}\right]\\
  &=e^R\E\left[e^{-Y_1/2}\right]^k\\ &\le h(R,\epsilon)^k
  e^R\end{split}\]
\noindent
completing the proof.\end{proof}
\noindent{\em Second part of proof of Theorem \ref{maintheorem2}.}
By the estimates in Proposition \ref{critboundlemma2} and Lemma
\ref{palbinlemma} we get that
\[
\sum_{l=2 k}^{\infty}\left(\lambda_c(R) \mu (S(0,2R))\right)^l
P_R(l)\le e^R\sum_{l=2k}^{\infty}K^l h(R,\epsilon)^{l-1}\]
\noindent
for any $\epsilon\in(0,1)$ and some constant $K\in(0,\infty)$.
Thus if we take $R$ big enough, the sum goes to $0$ as
$k\rightarrow\infty$. This is also the case if we replace
$\lambda_c$ with $t \lambda_c$ for some $t>1$, proving that
there are intensities above $\lambda_c$ for which there are
infinitely many unbounded connected components in the covered region
of $\Hypn$ for $R$ big enough. \qed

\bigskip\noindent {\bf Acknowledgements:} I want to thank Johan
Jonasson, my advisor, for introducing me to the problems dealt with in
this paper, and for all the help received while preparing the
manuscript. Thanks also to Olle H\"aggstr\"om for providing helpful
comments on parts of the manuscript. Finally I want to express my
thanks to Patrik Albin for giving permission to include the proof of Lemma
\ref{palbinlemma}.


\begin{thebibliography}{99}

\bibitem{albin} P. Albin, \textit{Private communication}.

\bibitem{alexander} K.S. Alexander, \textit{The RSW theorem for
    continuum percolation and the CLT for Euclidean minimal spanning
    trees}, Ann. Appl. Probab. {\bf 6} (1996), no. 2, 466-494.

\bibitem{athreya} K.B. Athreya and P.E. Ney, \textit{Branching
    Processes}, Springer Verlag, 1972.

\bibitem{itai2} I. Benjamini, R. Lyons, Y. Peres, and O. Schramm,
  \textit{Critical percolation on any nonamenable group has no
    infinite clusters}, Ann. Probab. {\bf 27} (1999), 1347-1356.

\bibitem{itai3} I. Benjamini, R. Lyons, Y. Peres, and O. Schramm,
  \textit{Group-invariant percolation on graphs}, Geom. Funct. Anal.
  {\bf 9} (1999), 29-66.

\bibitem{itai1} I. Benjamini and O. Schramm, \textit{Percolation in
    the hyperbolic plane}, J. Amer. Math. Soc. {\bf 14} (2001),
  487-507.

\bibitem{itai4} I. Benjamini and O. Schramm, \textit{Percolation
    beyond ${\mathbb Z}^d$, many questions and a few answers},
  Electronic Commun. Probab. {\bf 1} (1996), 71-82.

\bibitem{burtonkeane} R. M. Burton and M. Keane, \textit{Density and
    uniqueness in percolation}, Comm. Math. Phys. {\bf 121} (1989),
  501-505.

\bibitem{flavors} J.W. Cannon, W.J. Floyd, R. Kenyon and W.R. Parry.
  Hyperbolic geometry. In \emph{Flavors of geometry}, pp. 59-115,
  Cambridge University Press, 1997.

\bibitem{ErdMagObeTri0} A. Erd\'elyi, W. Magnus, F. Oberhettinger, and
  F.G.  Tricomi, {\it Higher Transcendental Functions, Vol.\
    I}, McGraw-Hill, 1953.

 \bibitem{ErdMagObeTri} A. Erd\'elyi, W. Magnus, F. Oberhettinger, and
   F.G. Tricomi, {\it Higher Transcendental Functions, Vol.\
     II}, McGraw-Hill, 1953.

\bibitem{grimmett} G. Grimmett, \textit{Percolation (2nd ed.)},
  Springer-Verlag, 1999.

\bibitem{ollehaegg} O. H\"aggstr\"om, \textit{Infinite clusters in
    dependent automorphism invariant percolation on trees}, Ann.
  Probab. {\bf 25} (1997), 1423-1436.

\bibitem{jonasson2} O. H\"aggstr\"om and J. Jonasson,
  \textit{Uniqueness and non-uniqueness in percolation theory},
  Probability Surveys {\bf 3} (2006), 289-344.

\bibitem{hall} P. Hall, \textit{On continuum percolation}, Ann.
  Probab. {\bf 13} (1985), 1250-1266.

\bibitem{jonasson} J. Jonasson, \textit{Hard-sphere percolation: Some
    positive answers in the hyperbolic plane and on the integer
    lattice}, Preprint, 2001.

\bibitem{paksmirnova} I. Pak and T. Smirnova-Nagnibeda, \textit{On
    non-uniqueness of percolation on nonamenable Cayley graphs}, C. R.
  Acad. Sci. Paris Sér. I Math {\bf 330} (2000), 495-500.

\bibitem{meester} R. Meester and R. Roy, \textit{Continuum
    Percolation}, Cambridge University Press, New York, 1996.

\bibitem{PruBryMar} A.P Prudnikov, Yu. A. Brychkov, and O.I. Marichev,
  {\it Integrals and Series. Volume 1: Elementary Functions,} Gordon
  and Breach Science Publishers, 1986.

\bibitem{ratcliffe} J.G. Ratcliffe, \textit{Foundations of hyperbolic
    manifolds}, Springer Verlag 2006.

\bibitem{sarkar} A. Sarkar, \textit{Co-existence of the occupied and
    vacant phase in Boolean models in three or more dimensions}, Adv.
  Appl. Prob. {\bf 29} (1997), 878-889.

\bibitem{schonmann} R.H. Schonmann, \textit{Stability of infinite
    clusters in supercritical percolation}, Probab. Th. Rel. Fields
  {\bf 113} (1999), 287-300.

\end{thebibliography}
\end{document}